\documentclass{article}

\usepackage{fullpage}
\usepackage{amsmath,amsfonts,amssymb,amsthm,enumerate,hyperref,cleveref,color,calc}

\DeclareMathOperator{\Tr}{Tr}
\DeclareMathOperator{\rank}{rank}
\DeclareMathOperator*{\Ex}{\mathbb{E}}
\DeclareMathOperator{\Var}{Var}
\DeclareMathOperator{\Kn}{Kn}
\providecommand{\symdiff}{\oplus}
\providecommand{\compl}{\overline}
\providecommand{\bfHD}{\mathbf{H}_{\mathcal{D} \times \mathcal{D}}}
\providecommand{\flt}{\tilde}
\providecommand{\bH}{\flt{\mathbf{H}}}
\providecommand{\bbfHD}{\bH_{\mathcal{D} \times \mathcal{D}}}

\newcommand{\remove}[1]{}

\newtheorem{theorem}{Theorem}

\title{Chvátal's conjecture: a proof from The Book}
\author{David Ellis, Yuval Filmus, Ehud Friedgut}
\date{September 2026}

\begin{document}

\maketitle

\begin{abstract}
Chvátal conjectured that every downset has a maximum-size intersecting family which is a star, that is, consists of all members of the downset containing a fixed element.

Recently, Chang, Liu and Liu gave a proof of this conjecture, which follows as a corollary of more general results such as Kleitman's conjecture and a version of Kahn's conjecture (which they also prove). We give a short, direct (spectral) proof of Chvátal's conjecture, and a proof of Kleitman's conjecture along the same lines. We also prove a strengthening of Chv\'atal's conjecture concerning the projection packing number. Finally, we propose two spectral Chv\'atal conjectures which are motivated by extensive numerical experiments.
\end{abstract}

\section{Introduction}
\label{sec:introduction}

Let $\mathcal{D}$ be a (finite) downset: a family of subsets of $[n]$, for some $n \in \mathbb{N}$, such that if $A \in \mathcal{D}$ and $B \subseteq A$ then $B \in \mathcal{D}$. An \emph{intersecting family} in $\mathcal{D}$ is a subset $\mathcal{I} \subseteq \mathcal{D}$ which is intersecting: if $A,B \in \mathcal{I}$ then $A \cap B \neq \emptyset$. What is the maximum possible size of an intersecting family in $\mathcal{D}$?

Chv\'atal~\cite{Chvatal74} conjectured that the maximum size of an intersecting family is always achieved by some star:

\begin{theorem} \label{thm:chvatal}
Let $\mathcal{D}$ be a (finite) downset. Then $\mathcal{D}$ has an intersecting family of maximum size which is a star.
\end{theorem}
This conjecture was one of the most famous open problems in extremal combinatorics, and despite half a century of effort, it was known only in (very) special cases; a list of more than 20 papers directly related to the conjecture appears at
\url{https://users.encs.concordia.ca/~chvatal/conjecture.html}.

Recently, Chang, Liu and Liu~\cite{CLL26} proved a correlation inequality conjectured by Friedgut, Kahn, Kalai and Keller~\cite{FKKK18} which implies Chv\'atal's conjecture. In fact, they also proved a stronger conjecture, due to Kleitman~\cite{Kleitman79}, that any antipodal intersecting family of subsets of $[n]$ flows down to a convex combination of stars, and a Boolean version of a beautiful conjecture of Kahn~\cite[Conjecture 3.6]{FKKK18}, which gives a prescribed convex combination of stars that satisfies the conclusion of Kleitman's conjecture. Shortly after, Keevash~\cite{Keevash26} used some of the methods of Chang, Liu and Liu to prove Kahn's conjecture in complete generality (i.e., for all increasing antipodal real-valued functions, not necessarily Boolean-valued).

The purpose of this note is threefold. First, we provide a ``Book Proof'' of Chv\'atal's conjecture, using ideas from the work of Chang, Liu and Liu. Second, we prove a new strengthening of Chv\'atal's conjecture, which bounds the projection packing number; the proof is a simple modification of the Book Proof. Finally, we report on empirical findings which suggest two other strengthenings which we were unable to resolve.

The same proof technique can also be used to simplify the proofs of the stronger results proved by Chang, Liu and Liu and by Keevash: Kleitman's conjecture and both versions of Kahn's conjecture. One can also formulate and prove analogs of these results in the setting of projection packings. In Section 3 we give the short proof of Kleitman's conjecture, along these lines. In order to keep this note short, we refrain from spelling out the details for the other generalizations, leaving them to the interested reader.

We note that it is usually inappropriate for mathematicians to refer to one of their own proofs as a ``Book Proof''; however, in this case we allow ourselves this luxury, due to the fact that the proof was not conceived by us, but rather coaxed by us out of an artificial medium; see our declaration at the end of the paper.

\section{Proof from The Book}
\label{sec:book-proof}

In this section we prove \Cref{thm:chvatal}. Let $\mathcal{D}$ be a downset of subsets of $[n]$. Throughout, $\compl{A}=[n]\setminus A$, and expectations are taken with respect to the uniform measure on $2^{[n]}$. We need to prove that if $\mathcal{I} \subseteq \mathcal{D}$ is intersecting then
\[
 |\mathcal{I}| \le \max_{i \in [n]} |\mathcal{S}_i|, \text{ where } \mathcal{S}_i = \{ A \in \mathcal{D} : i \in A \}.
\]
The proof will involve the \emph{upset generated by} $\mathcal{I}$:
\[
 \mathcal{I}^\uparrow := \{A \subseteq [n] : A \supseteq B \text{ for some } B \in \mathcal{I}\}.
\]

The family $\mathcal{I}^\uparrow$ is clearly intersecting. We may therefore define a function $h\colon 2^{[n]} \to \mathbb{R}$ as follows:
\[
 h(T) =
 \begin{cases}
     +1 & \text{if } T \in \mathcal{I}^\uparrow, \\
     -1 & \text{if } \compl{T} \in \mathcal{I}^\uparrow, \\
     \phantom{+}0 & \text{otherwise}.
 \end{cases}
\]
Clearly $\Ex[h] = (|\mathcal{I}^\uparrow| - |\mathcal{I}^\uparrow|)/2^n = 0$.

As usual (see e.g.\ \cite{ODonnell}), for any function $f\colon 2^{[n]} \to \mathbb{R}$ and any subset $S \subseteq [n]$, we let $\hat{f}(S) = \frac{1}{2^n}\sum_{T \subseteq [n]} f(T)\chi_S(T)$, where $\chi_S(T) = (-1)^{|S \cap T|}$; then $f$ has the Fourier expansion 
\[ f = \sum_{S \subseteq [n]}\hat{f}(S)\chi_S. \]
As usual, we write $L^2(2^{[n]})$ for the space of real-valued functions on $2^{[n]}$. Define a linear operator $\mathbf{H}\colon L^2(2^{[n]}) \to L^2(2^{[n]})$ as follows:
\[
 \widehat{\mathbf{H}f}(S) = h(S) \hat{f}(S)\quad \forall S \subseteq [n].
\]
Below, we also think of $\mathbf{H}$ as a $2^n \times 2^n$ matrix, identifying functions $2^{[n]} \to \mathbb{R}$ with vectors of length $2^n$.

The theorem will follow by combining two inequalities:
\[
 2|\mathcal{I}| \stackrel{\text{(LB)}}{\le} \sum_{A,B \in \mathcal{D}} \mathbf{H}(A,B)^2 \stackrel{\text{(UB)}}{\le} 2\max_{i \in [n]} |\mathcal{S}_i|.
\]

\paragraph{Proof of (UB):} Multiplication in the Fourier domain corresponds to convolution in the physical domain, and so
\[
 (\mathbf{H} f)(A) = \sum_B \hat{h}(A \symdiff B) f(B),
\]
where $\symdiff$ is symmetric difference. Therefore $\mathbf{H}(A,B) = \hat{h}(A \symdiff B)$,
and consequently
\[
 \sum_{A,B \in \mathcal{D}} \mathbf{H}(A,B)^2 = \sum_{A,B \in \mathcal{D}} \hat{h}(A \symdiff B)^2.
\]

Recall that $\hat{h}(\emptyset) = \Ex[h] = 0$. Therefore
\[
 \sum_{A,B \in \mathcal{D}} \hat{h}(A \symdiff B)^2 =
 \sum_{T \neq \emptyset} |\mathcal{D} \cap (\mathcal{D} \symdiff T)| \cdot \hat{h}(T)^2 \le 
 \max_{T \neq \emptyset} |\mathcal{D} \cap (\mathcal{D} \symdiff T)| \cdot \sum_T \hat{h}(T)^2 
 \stackrel{(\ast)}\le \max_{T \neq \emptyset} |\mathcal{D} \cap (\mathcal{D} \symdiff T)|,
\]
where $(\ast)$ follows from Parseval, $\hat{h}(\emptyset)=0$ and $\|h\|_2^2  \le 1$; here $\mathcal{D} \symdiff T := \{A \symdiff T : A \in \mathcal{D}\}$.

We now relate $|\mathcal{D} \cap (\mathcal{D} \symdiff T)|$ to the size of stars, for non-empty $T$. Let $T \neq \emptyset$ and choose any $i \in T$. Consider an arbitrary $A \in \mathcal{D} \cap (\mathcal{D} \symdiff T)$:
\begin{itemize}
\item If $A$ contains $i$ then $A \in \mathcal{S}_i$, hence there are at most $|\mathcal{S}_i|$ such sets $A$.
\item If $A$ doesn't contain $i$ then $A \symdiff T$ contains $i$ and so $A \symdiff T \in \mathcal{S}_i$, hence there are at most $|\mathcal{S}_i|$ such sets $A$.
\end{itemize}
Therefore $|\mathcal{D} \cap (\mathcal{D} \symdiff T)| \le 2|\mathcal{S}_i|$, completing the proof of (UB).

\paragraph{Proof of (LB):}

In order to prove (LB), we reinterpret the sum as a trace. Let $\bfHD$ be the restriction of $\mathbf{H}$ (viewed as a matrix) to the rows and columns indexed by $\mathcal{D}$. Then
\begin{equation}\label{eq:sum}
 \sum_{A,B \in \mathcal{D}} \mathbf{H}(A,B)^2 = \Tr\, ((\bfHD)^2) = \sum_{i=1}^{|\mathcal{D}|}\lambda_i^2,
\end{equation}
where $\lambda_1,\ldots,\lambda_{|\mathcal{D}|}$ are the eigenvalues of $\bfHD$, listed with their multiplicities. 

We will show that $\bfHD$ has at least $|\mathcal{I}|$ eigenvalues equal to $+1$ and at least $|\mathcal{I}|$ eigenvalues equal to $-1$. Consequently, the sum in \eqref{eq:sum} is at least $2|\mathcal{I}|$, completing the proof of (LB).

The matrix interpretation above immediately implies that if $g\colon 2^{[n]} \to \mathbb{R}$ is supported on $\mathcal{D}$, and $g$ is an eigenfunction of $\mathbf{H}$ with eigenvalue $\lambda$, then $g|_{\mathcal{D}}$ is an eigenfunction of $\bfHD$ with eigenvalue $\lambda$.

For $A \in \mathcal{I}$, let $p_A\colon 2^{[n]} \to \mathbb{R}$ be the characteristic function of the family of subsets of $A$, which is clearly supported on $\mathcal{D}$. As a function on $\{0,1\}^n$, $p_A$ corresponds to the statement that all coordinates in $\compl{A}$ are zero, and so its Fourier expansion is supported on subsets of $\compl{A}$. The definition of $\mathbf{H}$ now immediately implies that $\mathbf{H} p_A = -p_A$. It is easy to check that the functions $\{ p_A : A \in \mathcal{I} \}$ are linearly independent (if one views them as vectors, and one orders the entries in a manner consistent with containment, then the $|\mathcal{I}| \times |\mathcal{I}|$ matrix formed by the vectors $(p_A|_{\mathcal{I}})_{A \in \mathcal{I}}$ is upper triangular with non-zero entries on the diagonal). Since each $p_A$ is supported on $\mathcal{D}$, the restrictions $p_A|_{\mathcal{D}}$ are also linearly independent, and by the argument above they are eigenvectors of $\bfHD$ with eigenvalue $-1$, so we obtain $|\mathcal{I}|$ eigenvalues of $\bfHD$ equal to $-1$.

We obtain the $+1$ eigenvalues by considering the functions $q_A := \chi_{[n]} p_A$, for each $A \in \mathcal{I}$. The Fourier expansion of $q_A$ is supported on supersets of $A$, and so looking again at the definition of $\mathbf{H}$, we see that $\mathbf{H} q_A = q_A$. Each $q_A$ is supported on $\mathcal{D}$, and since $\chi_{[n]}$ is nowhere zero, the same ``upper triangular'' argument as for the $p_A$'s shows that the $q_A$'s (and therefore their restrictions to $\mathcal{D}$) are linearly independent. This gives $|\mathcal{I}|$ eigenvalues of $\bfHD$ equal to $+1$, completing the proof of (LB), and proving the theorem.

\section{A spectral proof of Kleitman's conjecture}
In this section, we use the spectral method of the ``Book Proof'' in the previous section to give a short proof of Kleitman's conjecture, in the following form. (Throughout this section, we write $\mathcal{S}_i:=\{A\subseteq[n]:i\in A\}$.)

\begin{theorem}\label{thm:kleitman}
Let $\mathcal{U}\subseteq2^{[n]}$ be a maximal intersecting family (meaning, no strict superset of $\mathcal{U}$ is intersecting). Then there exist real numbers $c_1,\ldots,c_n\ge0$, with $\sum_i c_i=1$, such that $1_{\mathcal{U}}$ flows down to $\sum_{i=1}^{n}c_i 1_{\mathcal{S}_i}$, or equivalently (by the max-flow min-cut theorem), for every downset $\mathcal{D} \subseteq 2^{[n]}$ we have 
\[
 |\mathcal{U} \cap \mathcal{D}| \le\sum_{i=1}^n c_i|\mathcal{S}_i \cap \mathcal{D}|.
\]
Furthermore, the coefficients $c_i$ can be defined as in (the Boolean case of) Kahn's conjecture.
\end{theorem}

We remind the reader that if $f,g\colon 2^{[n]} \to \mathbb{R}_{\geq 0}$, we say that $f$ \emph{flows down to} $g$ if there exists $v\colon 2^{[n]} \times 2^{[n]} \to \mathbb{R}_{\geq 0}$ such that for any $A \subseteq [n]$ we have $\sum_{B \subseteq [n]}v(A,B) = f(A)$, for any $B \subseteq [n]$ we have $\sum_{A \subseteq [n]}v(A,B)=g(B)$, and for any $B \not\subseteq A$, we have $v(A,B)=0$. (Cf.\ Definition 3.1 in \cite{FKKK18}, which defines an upward flow.)

We now describe Kahn's (conjectured) coefficients $c_i$. Let $\mathcal{U} \subseteq 2^{[n]}$ be a maximal intersecting family. Then $\mathcal{U}$ is an upset and contains exactly one member of every pair $\{A,\compl{A}\}$. Thus, setting $h=2\mathbf1_{\mathcal{U}}-1$, we have
\[
 h(\compl A)=-h(A)\quad \forall A \subseteq [n],\qquad
 \widehat h(\emptyset)=0,\qquad
 \sum_{T\ne\emptyset}\widehat h(T)^2=\Ex[h^2]=1.
\]
Following Kahn~\cite[Conjecture 3.4]{FKKK18}, for $i \in [n]$ define
\begin{equation*}
 c_i=\sum_{\substack{T\ne\emptyset\\ \max(T)=i}}\widehat h(T)^2.
\end{equation*}
These coefficients are clearly nonnegative and sum to one.

Here is an explanation (due to Kahn, personal communication) that gives a heuristic justification for these coefficients. For $i=0,1,\ldots,n$, define
\[
h_i(x)= \Ex[h(x) \mid x_1,x_2,\ldots,x_i].
\]
The sequence of $h_i$'s provides finer and finer approximations of $h$, with nondecreasing variance. The coefficient $c_i$ measures how much more detailed this approximation becomes when additionally conditioning on the $i$th coordinate, in the sense that
\[
c_i = \Var(h_i) - \Var(h_{i-1}).
\]

\begin{proof}[Proof of \Cref{thm:kleitman}]
Let $\mathcal{U} \subseteq 2^{[n]}$ be a maximal intersecting family, and let $h$ and the $c_i$ be defined as above. Writing $\mathcal{I} = \mathcal{U} \cap \mathcal{D}$, we have (exactly as in the proof of \Cref{thm:chvatal})

\begin{equation}\label{eq:lower}
 2|\mathcal{U}\cap\mathcal{D}| = 2|\mathcal{I}|
 \le \sum_{A,B\in\mathcal{D}}\widehat h(A\symdiff B)^2.
\end{equation}

As in the proof of \Cref{thm:chvatal}, if $i \in T$ then $|\mathcal{D} \cap (\mathcal{D} \symdiff T)| \le 2|\mathcal{D} \cap \mathcal{S}_i|$. Applying this inequality with $i = \max(T)$, we obtain
\begin{align*}
 \sum_{A,B\in\mathcal{D}}\widehat h(A\symdiff B)^2
 &=\sum_{T\ne\emptyset}|\mathcal{D}\cap(\mathcal{D}\symdiff T)|\widehat h(T)^2\\
&\le2\sum_{T\ne\emptyset}|\mathcal{D}\cap\mathcal{S}_{\max(T)}|\widehat h(T)^2\\
 &=2\sum_i c_i|\mathcal{D}\cap\mathcal{S}_i|.
\end{align*}
Combining this with \eqref{eq:lower} proves the theorem.
\end{proof}

Our proof shows that the ``max'' in Kahn's choice of coefficients is somewhat arbitrary. Indeed, for any ``choice'' function $\tau\colon 2^{[n]} \setminus \{\emptyset\} \to [n]$ such that $\tau(T) \in T$ for all $T$, the above proof shows that the following coefficients would work:
\[
 c_i^\tau = \sum_{T \,:\, \tau(T)=i} \hat{h}(T)^2.
\]
In fact, though, this apparently more general result follows from \Cref{thm:kleitman}, by applying (for each downset $\mathcal{D}$) the inequality in \Cref{thm:kleitman} where we reorder the coordinates so that $|\mathcal{D} \cap \mathcal{S}_1| \ge \cdots \ge |\mathcal{D} \cap \mathcal{S}_n|$.

\smallskip

As observed in \cite{FKKK18}, Kleitman's conjecture implies a $p$-biased analogue of Chv\'atal's conjecture: namely (defining, as usual, the $p$-biased measure $\mu_p$ on $\{0,1\}^{[n]}$ by $\mu_p(S) = p^{|S|}(1-p)^{n-|S|}$), if $\mathcal{D} \subseteq \{0,1\}^{[n]}$ is a downset, $\mathcal{I} \subseteq \mathcal{D}$ is intersecting, and $0 \leq p \leq 1/2$, then there exists $i \in [n]$ such that $\mu_p(\mathcal{I}) \leq \mu_p(\mathcal{S}_i \cap \mathcal{D})$. Indeed, Kleitman's conjecture implies that if $\mathcal{U} \subseteq 2^{[n]}$ is a maximal intersecting family, then there exist $c_1,\ldots,c_n \geq 0$ with $\sum_{i=1}^nc_i=1$ such that for any monotone decreasing function $h\colon \{0,1\}^{[n]} \to \mathbb{R}_{\geq 0}$, we have \[ \Ex_{\mu_{1/2}}[fh] \leq \Ex_{\mu_{1/2}}[gh], \]
where $f := 1_{\mathcal{U}}$ and $g := \sum_{i=1}^{n}c_i1_{\mathcal{D}_i}$. Now let $\mathcal{D} \subseteq 2^{[n]}$ be any downset. Taking $h(S) = 2^n p^{|S|}(1-p)^{n-|S|}1_{S \in \mathcal{D}}$ in the above inequality (and noting that $h$ is monotone decreasing for $p \leq 1/2$) yields 
\[ \mu_p(\mathcal{U} \cap \mathcal{D}) = \Ex_{\mu_{p}}[f1_{\mathcal{D}}] = \Ex_{\mu_{1/2}}[fh] \leq \Ex_{\mu_{1/2}}[gh] = \Ex_{\mu_{p}}[g1_{\mathcal{D}}] = \sum_{i=1}^{n}c_i \mu_p(\mathcal{S}_i \cap \mathcal{D}), \]
where $f$ and $g$ are as above, and therefore there exists $i \in [n]$ such that $\mu_p(\mathcal{U} \cap \mathcal{D}) \leq \mu_p(\mathcal{S}_i \cap \mathcal{D})$, as required.

\section{Projection packing bound}
\label{sec:projection-packing-bound}

In this section, for a downset $\mathcal{D}$, we return to the convention $\mathcal{S}_i=\{A\in\mathcal{D}:i\in A\}$.

A \emph{projection packing} is a generalization of an independent set in a graph. There are two equivalent ways to describe this notion: using vector spaces, or using orthogonal projections. We start with the former.

A \emph{$d$-dimensional projection packing} for a finite graph $G = (X,E)$ assigns to each vertex $x \in X$ a subspace $V_x$ of $\mathbb{C}^d$, such that $V_x \perp V_y$ whenever $xy \in E$. The \emph{value} of such a packing is defined to be 
\[
 \frac{1}{d} \sum_{x \in X} \dim V_x.
\]

Letting $P_x$ denote the orthogonal projection onto $V_x$, we see that a $d$-dimensional projection packing of $G$ corresponds exactly to assigning to each $x \in X$ an orthogonal projection $P_x\colon \mathbb{C}^d \to \mathbb{C}^d$, subject to the condition that $P_x P_y = 0$ whenever $xy \in E$. The value is then equal to
\[
 \frac{1}{d} \sum_{x \in X} \rank P_x.
\]

The \emph{projection packing number} of $G$ is the supremum of these values over all $d$ and all $d$-dimensional projection packings.

A $1$-dimensional projection packing for $G$ is the same as an independent set in $G$, and the value of the projection packing is equal to the cardinality of the independent set. This shows that the projection packing number of a graph is an upper bound on its independence number.
(On the other hand, as we explain in Section \ref{sec:conjectures}, the projection packing number is bounded from above by the theta bound and the inertia bound.)

For a downset $\mathcal{D}$, consider the Kneser graph $\Kn(\mathcal{D})$ on $\mathcal{D}$, i.e., the graph with vertex-set $\mathcal{D}$ where two sets are joined by an edge iff they are disjoint; there is a loop on the vertex $\emptyset$. An independent set in $\Kn(\mathcal{D})$ is the same as an intersecting family in $\mathcal{D}$.

\begin{theorem} \label{thm:packing}
Let $\mathcal{D}$ be a (finite) downset. For every positive integer $d$, the largest value of a $d$-dimensional projection packing for $\Kn(\mathcal{D})$ is the size of the largest star.
\end{theorem}

Another related parameter is the \emph{quantum independence number} of a graph $G=(X,E)$. A \emph{quantum independent set of value $k$} in $G$ consists of $k$ orthogonal decompositions $(V_{1,x})_{x \in X},\dots,(V_{k,x})_{x \in X}$ of $\mathbb{C}^d$ for some $d \in \mathbb{N}$ (allowing zero subspaces), such that $V_{1,x},\dots,V_{k,x}$ are mutually orthogonal for all $x \in X$, and $V_{i,x} \perp V_{j,y}$ for all $i,j \in [k]$ whenever $xy \in E$. A quantum independent set of value $k$ immediately yields a projection packing of value $k$, by taking $V_x := \bigoplus_{i=1}^k V_{i,x}$ for all $x \in X$. The quantum independence number of $G$ is the largest $k$ for which there exists a quantum independent set of value $k$. Note that an independent set of size $k$ in $G$ gives rise to a quantum independent set of value $k$. \Cref{thm:packing} thus shows that the maximum size of an intersecting family, the quantum independence number of $\Kn(\mathcal{D})$, and the projection packing number of $\Kn(\mathcal{D})$ all coincide with the maximum size of a star.

We prove \Cref{thm:packing} in the rest of this section, closely following the proof in \Cref{sec:book-proof}.

\medskip

Let $(U_A)_{A \in \mathcal{D}}$ be a projection packing for $\Kn(\mathcal{D})$. We construct the corresponding ``upset'' (cf.~\Cref{sec:book-proof}) by defining, for all $A \subseteq [n]$,
\[
 V_A = \operatorname{span}(U_B : B \subseteq A, B \in \mathcal{D}).
\]
If $A \cap B = \emptyset$ then $V_A \perp V_B$. Moreover, if $A \in \mathcal{D}$ then $\dim V_A \ge \dim U_A$. Therefore $(V_A)_{A \in \mathcal{D}}$ is a projection packing whose value is at least as large as the value of $(U_A)_{A \in \mathcal{D}}$, and so it suffices to bound the value of $(V_A)_{A \in \mathcal{D}}$. The new projection packing has the useful property that $V_A \subseteq V_B$ whenever $A \subseteq B$. For every $A \subseteq [n]$, let $P_A$ be the orthogonal projection onto $V_A$.

The real-valued function $h$ is replaced by the $\mathbb{C}^{d \times d}$-valued function
\[
 h(T) = P_T - P_{\compl{T}}\quad \forall T \subseteq [n].
\]
When $d = 1$, we have $P_T \in \{0,1\}$, and so this agrees with the definition in \Cref{sec:book-proof}. By definition, $h(T)$ is Hermitian.

In order to define the linear operator $\mathbf{H}$, we first extend the Fourier transform to functions $f\colon 2^{[n]} \to \mathbb{C}^d$ by using the standard defining formula: $\hat{f}(S) = \frac{1}{2^n} \sum_{T \subseteq [n]} f(T) \chi_S(T) \in \mathbb{C}^d$. We now define the operator $\mathbf{H}$ acting on such functions using the same formula as in \Cref{sec:book-proof}:
\[
 \widehat{\mathbf{H}f}(S) = h(S) \hat{f}(S)\quad \forall S \subseteq [n].
\]

The theorem will follow by combining two inequalities:
\[
 2\sum_{A \in \mathcal{D}} \dim V_A \stackrel{\text{(LB)}}{\le} \sum_{A,B \in \mathcal{D}} \Tr \mathbf{H}(A,B)^2 \stackrel{\text{(UB)}}{\le} 2d\max_{i \in [n]} |\mathcal{S}_i|.
\]
Compared with \Cref{sec:book-proof}, we replace squared matrix entries by traces of squared blocks, replace $|\mathcal{I}|$ by $\sum_{A\in\mathcal{D}}\dim V_A$, and gain a factor of $d$ in the upper bound.

\paragraph{Proof of (UB):} Multiplication in the Fourier domain corresponds to convolution in the physical domain:
\[
 (\mathbf{H}f)(A) = \sum_B \hat{h}(A \symdiff B) f(B),
\]
where the $\mathbb{C}^{d \times d}$-valued Fourier expansion of $h$ is defined using the standard defining formula as above. Since $h(T)$ is Hermitian for every $T$, the Fourier coefficients $\hat{h}(S)$ are Hermitian for every $S$, and so the matrix $\mathbf{H}$ is Hermitian (that is, $\mathbf{H}(A,B) = \mathbf{H}(B,A)^*$).

Clearly $\hat{h}(\emptyset) = \Ex[h] = (\sum_T P_T - \sum_T P_{\compl{T}})/2^n = 0$.
As in \Cref{sec:book-proof}, this implies that
\[
 \sum_{A,B \in \mathcal{D}} \Tr \mathbf{H}(A,B)^2 =
 \sum_{A,B \in \mathcal{D}} \Tr \hat{h}(A \symdiff B)^2 \le 
 \max_{T \neq \emptyset} |\mathcal{D} \cap (\mathcal{D} \symdiff T)| \cdot \sum_T \Tr \hat{h}(T)^2.
\]

We showed in \Cref{sec:book-proof} that $\max_{T \neq \emptyset} |\mathcal{D} \cap (\mathcal{D} \symdiff T)| \le 2\max_i |\mathcal{S}_i|$. We complete the proof of (UB) by bounding $\sum_T \Tr \hat{h}(T)^2$ using Parseval's identity:
\[
 \sum_T \Tr \hat{h}(T)^2 = \Ex_T[\Tr h(T)^2] = \Ex_T[\Tr(P_T - P_{\compl{T}})^2] = \Ex_T[\Tr P_T + \Tr P_{\compl{T}}] \le d,
\]
using $P_T P_{\compl{T}} = 0$ and $\Tr P_T + \Tr P_{\compl{T}} = \dim V_T + \dim V_{\compl{T}} \le d$.

\paragraph{Proof of (LB):} Following the approach of \Cref{sec:book-proof}, we interpret the sum $\sum_{A,B \in \mathcal{D}} \Tr \mathbf{H}(A,B)^2$ as a single trace. To this end, recall that for a Hermitian matrix $\mathbf{M}$ (with complex-valued entries), the quantity $\Tr \mathbf{M}^2$ is the sum of the squared absolute values of all entries. Hence if $\bbfHD$ is the flattened $(|\mathcal{D}| \cdot d) \times (|\mathcal{D}| \cdot d)$ matrix (with complex-valued entries) corresponding to $\bfHD$, then
\[
 \sum_{A,B \in \mathcal{D}} \Tr \mathbf{H}(A,B)^2 = \Tr\, ((\bbfHD)^2) = \sum_{i=1}^{|\mathcal{D}| \cdot d} \lambda_i^2,
\]
where $\lambda_1,\dots,\lambda_{|\mathcal{D}|\cdot d}$ are the (real) eigenvalues of $\bbfHD$, listed with their multiplicities. 

Following the approach of \Cref{sec:book-proof}, we construct eigenvectors of $\bbfHD$ by constructing eigenvectors supported on $\mathcal{D} \times [d]$ of the flattened $(2^n \cdot d) \times (2^n \cdot d)$ complex-valued matrix $\bH$. To this end, it will be convenient to switch back to the original matrix $\mathbf{H}$: an eigenvector $\flt{f} \in \mathbb{C}^{2^n \cdot d}$ of $\bH$ corresponds to an eigenfunction $f\colon 2^{[n]} \to \mathbb{C}^d$ of $\mathbf{H}$. The eigenvector $\flt{f}$ is supported on $\mathcal{D} \times [d]$ iff the eigenfunction $f$ is supported on $\mathcal{D}$. We are now ready to construct the required eigenfunctions.

For every $A \in \mathcal{D}$, define the functions $p_A,q_A$ as in \Cref{sec:book-proof}. We will show that for each nonzero $v \in V_A$, the function $p_{A,v}\colon 2^{[n]} \to \mathbb{C}^d$ given by $p_{A,v}(T) := p_A(T) v$ is an eigenfunction of $\mathbf{H}$ with eigenvalue $-1$, and the function $q_{A,v}\colon 2^{[n]} \to \mathbb{C}^d$ given by $q_{A,v}(T) := q_A(T) v$ is an eigenfunction of $\mathbf{H}$ with eigenvalue $+1$. Since both are supported on $\mathcal{D}$, the restrictions of the corresponding vectors $\flt{p}_{A,v},\flt{q}_{A,v}$ to $\mathcal{D} \times [d]$ are eigenvectors of $\bbfHD$.

For every $A \in \mathcal{D}$, the subspace $\{ p_{A,v} : v \in V_A\}$ has dimension $\dim V_A$. Moreover, the sum of these subspaces is direct since the functions $p_A$ are linearly independent, and so the eigenspace of $-1$ has dimension at least $\sum_{A \in \mathcal{D}} \dim V_A$. Replacing $p_{A,v}$ with $q_{A,v}$, we reach a similar conclusion for the eigenspace of $+1$.

To prove that $q_{A,v}$ is an eigenfunction, observe that its Fourier expansion is supported on supersets $T \supseteq A$. Therefore
\[
 \widehat{\mathbf{H} q_{A,v}}(T) = (P_T - P_{\compl{T}}) \widehat{q_{A,v}}(T) = \widehat{q_A}(T) (P_T - P_{\compl{T}}) v = \widehat{q_A}(T) v = \widehat{q_{A,v}}(T).
\]
Here we use $v \in V_A \subseteq V_T$ and so $P_T v = v$, and $V_A \perp V_{\compl{T}}$ (since $A \cap \compl{T} = \emptyset$) and so $P_{\compl{T}} v = 0$.

Similarly, to prove that $p_{A,v}$ is an eigenfunction, observe that its Fourier expansion is supported on sets $T$ such that $\compl{T} \supseteq A$. Hence the calculation above shows that $\widehat{\mathbf{H} p_{A,v}}(T) = -\widehat{p_{A,v}}(T)$.

\section{Two spectral conjectures}
\label{sec:conjectures}

Although it is a well-established method in combinatorics and related fields to consider the trace of the square of a symmetric matrix (equivalently, the sum of the squares of its eigenvalues), and to bound this using eigenvalue multiplicities, we feel that some elements of the proof we have given of \Cref{thm:chvatal} are ad hoc, without a clear precedent in the literature (prior to the work of Chang, Liu and Liu) --- for example, the choice of $\mathbf{H}$ --- and we feel the same about the proof of Chang, Liu and Liu. It is natural to ask whether Chv\'atal's conjecture can be proved using more ``standard'' techniques, of the type used to bound the maximum size of an independent set in a graph, applied to the Kneser graph $\Kn(\mathcal{D})$ as defined in \Cref{sec:projection-packing-bound}. Two of them come to mind:

\begin{description}
\item[(Weighted) Hoffman bound:] Let $G = (V,E)$ be a finite graph,\footnote{For technical reasons, we allow vertices of $G$ to have loops (the vertex $\emptyset$ of $\Kn(\mathcal{D})$ has a loop). For the definition of the weighted Hoffman bound to make sense, we need at least one vertex without a loop. When discussing the Lovász theta number and fractional independence number later on, we delete vertices with loops.} and suppose $\mathbf{M}$ is a real symmetric $V \times V$ matrix such that $\mathbf{M}(x,y) = 0$ whenever $xy \notin E(G)$, and $\mathbf{M} \mathbf{1} = \mathbf{1}$, where $\mathbf{1}$ is the all-$1$ vector. Then the maximum size of an independent set in $G$ is bounded from above by $\frac{-\lambda_{\min}(\mathbf{M})}{1-\lambda_{\min}(\mathbf{M})} |V(G)|$, where $\lambda_{\min}(\mathbf{M})$ is the minimal eigenvalue of $\mathbf{M}$. 

\item[The inertia bound:] Let $G = (V,E)$ be a finite graph, and suppose $\mathbf{M}$ is a real symmetric $V \times V$ matrix such that $\mathbf{M}(x,y) = 0$ whenever $xy \notin E(G)$. Then the maximum size of an independent set in $G$ is bounded from above by the number of nonnegative eigenvalues of $\mathbf{M}$. 
\end{description}

For any graph, the projection packing number considered in \Cref{sec:projection-packing-bound} is bounded from above by the inertia bound and by the Lov\'asz theta bound; see~\cite{Roberson13,WEA22}. The theta bound coincides with the weighted Hoffman bound in the case where there exists a vertex with a loop that is joined to every other vertex (as is the case for the vertex $\emptyset$ in $\Kn(\mathcal{D})$), an observation appearing implicitly in~\cite[Theorems 3 and 9]{Lovasz79}. The weighted Hoffman bound and the inertia bound are incomparable; see~\cite{ihringer,kwan}.

Hoping to prove Chv\'atal's conjecture by somehow using the inertia bound or  the weighted Hoffman bound (specialized to $\Kn(\mathcal{D})$), we checked whether they were tight on thousands of downsets $\mathcal{D}$, both structured and random. Amazingly, both bounds worked on every family that we tested! We therefore make the following conjectures:

\medskip

\textbf{Conjecture H:} The weighted Hoffman bound is tight for any downset.

\medskip

\textbf{Conjecture \makebox[\widthof{H}]{I}:} The inertia bound is tight for any downset.

\medskip

Unfortunately, the new proofs of \Cref{thm:chvatal} do not seem to produce a matrix certifying a tight bound in either case.

We also considered a strengthening of Conjecture~H obtained by replacing the weighted Hoffman bound by the fractional independence number, which is an upper bound on the Lov\'asz theta number (and hence on the optimal weighted Hoffman bound). The fractional independence number is defined for a graph $G = (V,E)$ as the maximal value of $\sum_{v \in V} x_v$ over all nonnegative $x_v$ satisfying the constraints $\sum_{v \in C} x_v \le 1$ for every clique $C \subseteq V$. It turns out that this conjecture is false already for $\mathcal{D} = \{A \subseteq [7] : |A| \leq 3\}$, where the maximum intersecting family has size $22$, while the fractional independence number is $22.75$.

\paragraph{AI disclosure and chronology of our work on this project} 
We formulated Conjecture H circa 2010 (unpublished), following some numerical experimentation. Recently, we began attempting to prove it with the assistance of ChatGPT-6~Astra and Claude~Fable~5.1, but to no avail. After a couple of weeks of unsuccessful attempts, we came up with Conjecture I, and suggested it to the bots, as a complementary project. Our confidence in Conjecture I at that point was quite low, but to our astonishment it turned out to hold on thousands of examples. Here too, ChatGPT and Claude were unable to supply a proof. (From this point on we stopped using Claude, and used only ChatGPT.)

Our next attempt to make a breakthrough was via the quantum independence number and the projection packing number, both of which are bounded below by the independence number and bounded above by each of the inertia bound (of Conjecture~I) and the Lov\'asz theta bound (of Conjecture~H). Here too, even on this weaker conjecture, ChatGPT was unable to help us. This did not come as a surprise to us, since, as far as we know, there are no combinatorial problems where these parameters give an otherwise unknown upper bound on the independence number of a graph (not counting the example in this paper).

The crucial turning point came when Chang, Liu and Liu posted their proof.
We fed their manuscript into ChatGPT-6~Astra, and asked it to try again to prove Conjectures H and I, using these new insights. It didn't succeed, but to our astonishment it notified us that it now had a proof that the projection packing number is tight for every downset. After reading this proof we felt that some essential points were unclear, so we coaxed ChatGPT to strip this down to the spectral bare bones, and it then supplied us with the ``Book Proof'' presented in this paper.

We prepared the final paper completely on our own, using ChatGPT only for proofreading at this stage.


\bibliographystyle{alphaurl}
\bibliography{biblio}

\end{document}